\documentclass[11pt,reqno]{amsart}

\usepackage[utf8]{inputenc}
\usepackage[T1]{fontenc}
\usepackage{lmodern}
\usepackage[a4paper,margin=1.15in]{geometry}
\usepackage{amsmath,amssymb,amsthm,mathtools}
\usepackage{enumitem}
\usepackage{xcolor}
\usepackage[colorlinks=true,linkcolor=blue!60!black,citecolor=blue!60!black,urlcolor=blue!60!black]{hyperref}

\newtheorem{theorem}{Theorem}

\title[Asymptotic Gauss--Lucas theorem for finite Blaschke products]{Asymptotic Gauss--Lucas theorem for finite Blaschke products}

\author{Aneesh Jatar}
\address{Department of Mathematics, The University of Hong Kong, Pokfulam Road, Hong Kong, China}
\email{aneesh.jatar1@gmail.com}

\date{}

\begin{document}

\begin{abstract}
The Gauss--Lucas theorem states that if all the zeros of a polynomial lie in a compact convex set $K\subset\mathbb{C}$ then all zeros of its derivative also lie in $K$. An asymptotic version of this result was proved by Totik in \cite{Totik2016}, namely, if all but $o(n)$ zeros of a sequence of polynomials $P_n$, $n=1,2,\dots$ of degrees $\deg(P_n)=n$ lie in a compact convex set $K\subset\mathbb{C}$ then all but $o(n)$ zeros of their derivatives $P_n'$, $n=1,2,\dots$ lie in any fixed neighborhood of $K$. The Gauss--Lucas theorem has an analogue for finite Blaschke products, first established by J.~L.~Walsh in \cite{Walsh1952}, which states that if the zeros of a finite Blaschke product lie in a compact, hyperbolically convex set $K\subset\mathbb{D}$ then all zeros of its derivative in the unit disk $\mathbb{D}$ also lie in $K$. In this paper, we establish an analogue of Totik's asymptotic version of the Gauss--Lucas theorem for finite Blaschke products.
\end{abstract}

\maketitle

\section{Introduction}

The Gauss--Lucas theorem states that if all zeros of a complex polynomial $P\in\mathbb{C}[z]$ lie in a compact convex set $K\subset\mathbb{C}$, then all zeros of its derivative also lie in $K$. In \cite{Totik2016} the below asymptotic version of the Gauss--Lucas theorem was proved.

\begin{theorem}[Totik \cite{Totik2016}; 2016]\label{thm:totik}
Let $P_n$, $n=1,2,3,\dots$ be a sequence of polynomials, with $\deg(P_n)=n$. Suppose $K\subset\mathbb{C}$ is a compact, convex set, such that at most $o(n)$ of the zeros of $P_n$ lie outside $K$. Define for any fixed $\epsilon>0$ the $\epsilon$-thickening of $K$ to be the set $K_\epsilon:=\{z:\operatorname{dist}(z,K)<\epsilon\}$. Then for any fixed $\epsilon>0$, at most $o(n)$ of the zeros of $P_n'$ lie outside $K_\epsilon$.
\end{theorem}

A finite Blaschke product of degree $n\in\mathbb{N}_{\ge 1}$ is a rational function of the form
\[
B(z)=\lambda\prod_{j=1}^{n}\frac{z-z_j}{1-\overline{z_j}z},
\]
where $z_j$ is a sequence of points in $\mathbb{D}$ and $|\lambda|=1$ is a unimodular constant.

The finite Blaschke products may be thought of as the hyperbolic analogue of the complex polynomials, in the sense that both classes of functions are exactly the surjective, finite proper holomorphic mappings of their respective domains. These domains are $\mathbb{C}$ (in the case of complex polynomials) and $\mathbb{D}$ (in the case of finite Blaschke products). This point of view was emphasized by Ng and Tsang in \cite{NgTsang2013}, where a dictionary was established between complex polynomials and finite Blaschke products. In particular, given a result which holds for complex polynomials, based on the dictionary in \cite{NgTsang2013}, one would expect a hyperbolic analogue to hold for finite Blaschke products.

This correspondence was shown to extend to the probabilistic setting in \cite{CheungNgYam2014}, where an analogue of the Pemantle--Rivin conjecture (\cite{PemantleRivin2011}, Conjecture 1.1) was proved for finite Blaschke products. In this paper we prove an analogue of Theorem~\ref{thm:totik} for finite Blaschke products.

\begin{theorem}\label{thm:main}
Let $B_n$, $n=1,2,3,\dots$ be a sequence of finite Blaschke products, with $\deg(B_n)=n$. Suppose $K\subset\mathbb{D}$ is a hyperbolically convex, compact set such that at most $o(n)$ of the zeros of $B_n$ lie outside $K$. Define for any fixed $\epsilon>0$, the $\epsilon$-thickening of $K$ to be $K_\epsilon:=\{z\in\mathbb{D}:\operatorname{dist}_{\mathbb{D}}(z,K)<\epsilon\}$, where $\operatorname{dist}_{\mathbb{D}}(z,K)$ is the hyperbolic distance from $z$ to $K$. Then for any fixed $\epsilon>0$, at most $o(n)$ of the zeros of $B_n'$ contained in $\mathbb{D}$ lie outside $K_\epsilon$.
\end{theorem}

Our proof of Theorem~\ref{thm:main} essentially follows the strategy of Totik's proof of Theorem~\ref{thm:totik} in \cite{Totik2016} (which proceeds via a contradiction argument), while making necessary adjustments to accommodate for the hyperbolic setting.

\section{Proof of Theorem \ref{thm:main}}

Central to our argument is the Gauss--Lucas theorem for finite Blaschke products, first established by Walsh in \cite{Walsh1952}. We record this result below. For a rational function $B$ we let $Z(B)$ be its set of zeros in the unit disk $\mathbb{D}$ (with multiplicities taken into account).

\begin{theorem}[Walsh \cite{Walsh1952}]\label{thm:walsh}
Let $B$ be a finite Blaschke product with $Z(B)\subset K$, where $K\subset\mathbb{D}$ is a compact, hyperbolically convex set. Then $Z(B')\subset K$.
\end{theorem}

\begin{proof}[Proof of Theorem \ref{thm:main}]
Suppose for the sake of a contradiction that the statement of Theorem~\ref{thm:main} were false. Then (after passing to an appropriate subsequence whose relabeling we suppress) there is a value $0<\alpha<1$ such that for sufficiently large $n$, there are at least $\alpha n$ zeros of $B_n'$ lying outside $K_\epsilon$. Over the course of the proof of Theorem~\ref{thm:main}, we will frequently pass to a further subsequence, whose relabeling will always be suppressed. After passing to an appropriate further subsequence, we may assume that for the empirical measures
\[
\mu_n:=\frac{1}{n}\sum_{z\in Z(B_n)}\delta_z, \qquad
\nu_n:=\frac{1}{n}\sum_{z\in Z(B_n')}\delta_z,
\]
we have
\[
\mu_n\rightarrow_{n\rightarrow\infty}\mu
\]
and
\[
\nu_n\rightarrow_{n\rightarrow\infty}\nu
\]
in the weak$^*$-topology of measures on the Riemann sphere $\overline{\mathbb{C}}$.

Then necessarily, $\nu(K_\epsilon)\le\liminf_{n}\nu_n(K_\epsilon)\le 1-\alpha<1$, while $\mu(K)=1$. We will show that in fact $\nu(K)=1$, which will establish our desired contradiction.

We may assume that
\[
B_n(z)=\prod_{w\in Z(B_n)}\frac{z-w}{1-\overline{w}z},
\]
so that
\[
\frac{1}{n}\log\frac{1}{|B_n(z)|}=U_{\mathbb{D}}^{\mu_n}(z).
\]
Here we let, for any Borel measure $\rho$ with compact support in $\mathbb{D}$,
\[
U_{\mathbb{D}}^{\rho}(z):=\int_{\mathbb{D}}g_{\mathbb{D}}(z,w)\,d\rho(w)
\]
denote the Green potential of $\rho$, where $g_{\mathbb{D}}(z,w)$ is the value at $z\in\mathbb{D}$ of the Green function of the unit disk, with pole at the point $w\in\mathbb{D}$. It is easy to see that, up to a harmonic function, the Green potential is the same as the logarithmic potential, in particular we still have the distributional equality
\[
-\frac{1}{2\pi}\Delta U_{\mathbb{D}}^{\rho}=\rho.
\]
For any Borel measure $\rho$, with compact support in $\mathbb{C}$, we define the logarithmic potential of $\rho$ to be
\[
U^{\rho}(z):=\int_{\mathbb{C}}\log\frac{1}{|z-w|}\,d\rho(w)\in L_{\mathrm{loc}}^{1}(\mathbb{C}).
\]

The remainder of this note will use some basic concepts of logarithmic potential theory. In particular we will make use of the notion of the Balayage or sweeping-out of a Borel measure onto a compact subset of the plane (see Section 2.4 of \cite{SaffTotik2024} for the definition of Balayage).

By the Poisson--Jensen formula (see Theorem 2.4.10, page 105 of Totik's book \cite{SaffTotik2024}), letting
\[
C_n(z):=\prod_{w\in Z(B_n')}\frac{z-w}{1-\overline{w}z},
\]
we have the equality:
\[
\frac{1}{n}\log\left|\frac{B_n'(z)}{n}\right|=\frac{1}{n}\log|C_n(z)|+\frac{1}{n}\int_{-\pi}^{\pi}\frac{1-|z|^2}{2\pi|z-\exp(i\theta)|^2}\log\left|\frac{B_n'(\exp(i\theta))}{n}\right|d\theta.
\]
We now show that the term
\[
\Delta_n(z):=\frac{1}{n}\int_{-\pi}^{\pi}\frac{1-|z|^2}{2\pi|z-\exp(i\theta)|^2}\log\left|\frac{B_n'(\exp(i\theta))}{n}\right|d\theta
\]
is $o(1)$ locally uniformly in $\mathbb{D}$, as $n\to\infty$, after passing to a further subsequence. By the hypothesis of Theorem~\ref{thm:main}, for any $\sup_{z\in K}|z|<r<1$ fixed, we have that at least $m_n$ of the zeros of $B_n$ remain inside $\mathbb{D}_r$, so that
\[
\frac{|B_n'(\exp(i\theta))|}{n}=\frac{1}{n}\sum_{w\in Z(B_n)}\frac{1-|w|^2}{|w-\exp(i\theta)|^2}\ge\frac{1}{n}\sum_{w\in Z(B_n)\cap\mathbb{D}_r}\frac{1-r^2}{(1+r)^2}=\frac{m_n}{n}\frac{1-r}{1+r},
\]
where $\frac{m_n}{n}\to 1$ as $n\to\infty$, hence
\[
\Delta_n(z)\ge\frac{1}{n}\log\left[\frac{m_n}{n}\frac{1-r}{1+r}\right],
\]
and $n\Delta_n(z)$ is a sequence of harmonic functions which are uniformly bounded below on the unit disk. On the other hand,
\[
n\Delta_n(0)=\int_{-\pi}^{\pi}\log\left|\frac{B_n'(\exp(i\theta))}{n}\right|\frac{d\theta}{2\pi}\le\int_{-\pi}^{\pi}\left|\frac{B_n'(\exp(i\theta))}{n}\right|\frac{d\theta}{2\pi}=\frac{1}{n}\int_{-\pi}^{\pi}\sum_{w\in Z(B_n)}\frac{1-|w|^2}{|w-\exp(i\theta)|^2}\frac{d\theta}{2\pi}=1.
\]
Since $n\Delta_n(0)$ is uniformly bounded below and above, by passing to a further subsequence, we may assume that $n\Delta_n(0)$ converges to a finite (real) value. Then by Harnack's theorem, which is applicable here as $n\Delta_n$ are uniformly bounded below, we may, after passing to a further subsequence, assume that $n\Delta_n$ converges locally uniformly on the unit disk to a harmonic function (finite-valued), which implies that $\Delta_n=o(1)$ holds locally uniformly in $\mathbb{D}$.

Now because
\[
U_{\mathbb{D}}^{\nu_n}(z)=-\frac{1}{n}\log|C_n(z)|,
\]
we obtain that
\[
\frac{1}{n}\log\left|\frac{B_n'(z)}{nB_n(z)}\right|=U_{\mathbb{D}}^{\mu_n}(z)-U_{\mathbb{D}}^{\nu_n}(z)+o(1).
\]
In other words, the left-hand side is up to a vanishing contribution, equal to the difference of Green potentials of the empirical measures of zeros and critical points.

The main idea for the rest of the proof is as follows:
\begin{enumerate}[label=\arabic*., leftmargin=*]
\item Show that in the limit (after passing to an appropriate further subsequence), the left-hand side
\[
\frac{1}{n}\log\left|\frac{B_n'(z)}{nB_n(z)}\right|
\]
of the above equality vanishes (Lebesgue) almost everywhere in $\mathbb{D}\setminus K$. This will follow from an approximation argument combined with Walsh's Gauss--Lucas theorem for finite Blaschke products (Theorem~\ref{thm:walsh}), details of this argument are provided on pages~\pageref{step1-start}--\pageref{step1-end}.

\item Letting $\widehat{\mu_n}$ and $\widehat{\nu_n}$ denote the Balayages of $\mu_n$ and $\nu_n$ onto $\mathbb{T}^1$ respectively, we may without loss of generality (after passing to an appropriate further subsequence), assume that $\widehat{\mu_n}\to\eta_1$ and $\widehat{\nu_n}\to\eta_2$ (in the weak$^*$-topology of Borel measures on the Riemann sphere $\overline{\mathbb{C}}$), for some probability measures $\eta_1,\eta_2$ supported on $\mathbb{T}^1$.

Since
\[
U_{\mathbb{D}}^{\mu_n}(z)=U^{\mu_n}(z)-U^{\widehat{\mu_n}}(z)
\]
and
\[
U_{\mathbb{D}}^{\nu_n}(z)=U^{\nu_n}(z)-U^{\widehat{\nu_n}}(z)
\]
(see Theorem X.1.1 of \cite{SaffTotik2024}, page 454) we may assume without loss of generality (again passing to a further subsequence) that
\[
U_{\mathbb{D}}^{\mu_n}\rightarrow_{n\rightarrow\infty}U^{\mu}-U^{\eta_1}
\]
and
\[
U_{\mathbb{D}}^{\nu_n}\rightarrow_{n\rightarrow\infty}U^{\nu}-U^{\eta_2}
\]
in $L_{\mathrm{loc}}^{1}(\mathbb{D})$ and almost everywhere in $\mathbb{D}$. Note that $U^{\eta_i}$, $i=1,2$ are both harmonic in $\mathbb{D}$, since $\eta_i$, $i=1,2$ are probability measures supported on $\mathbb{T}^1$.

\item By Steps 1 and 2, we conclude that the equality
\[
U^{\mu}=U^{\nu}+h,
\]
where $h=U^{\eta_1}-U^{\eta_2}$ is harmonic on all of $\mathbb{D}$, holds (Lebesgue) almost everywhere in $\mathbb{D}\setminus K$. Since $U^{\mu}$ is harmonic in $\mathbb{D}\setminus K$, its distributional Laplacian vanishes there, so that the distributional Laplacian of $U^{\nu}$ also vanishes there, and we conclude that $\nu$ has no mass in $\mathbb{D}\setminus K$. Using the fact that
\[
U^{\rho}(z_0)=\lim_{r\to 0}\frac{1}{\pi r^2}\int_{\Delta(z_0,r)}U^{\rho}(w)\,dA(w)
\]
holds at every point $z_0\in\mathbb{C}$ for any compactly supported (in the complex plane $\mathbb{C}$) Borel measure $\rho$, where $dA$ denotes the two-dimensional Lebesgue measure, the almost-everywhere equality $U^{\mu}=U^{\nu}+h$ in $\mathbb{D}\setminus K$ implies that in fact, $U^{\mu}(z)=U^{\nu}(z)+h(z)$ holds at every point $z\in\mathbb{D}\setminus K$.

\item Since $U^{\mu}(z)=U^{\nu}(z)+h(z)$ holds for every $z\in\mathbb{D}\setminus K$, and $\nu$ is compactly supported (in the complex plane) with no mass in $\mathbb{D}\setminus K$, letting $\gamma$ be any analytic Jordan curve in $\mathbb{D}\setminus K$ encircling $K$ (i.e.\ $K\subset\operatorname{int}(\gamma)$), the below Gauss' Theorem combined with the fact that $h=U^{\eta_1}-U^{\eta_2}$ is harmonic on all of $\mathbb{D}$, gives
\[
\nu(\operatorname{int}(\gamma))=\frac{1}{2\pi}\int_{\gamma}\frac{\partial U^{\nu}}{\partial\hat{n}}\,ds=\frac{1}{2\pi}\int_{\gamma}\frac{\partial(U^{\mu}-h)}{\partial\hat{n}}\,ds=\frac{1}{2\pi}\int_{\gamma}\frac{\partial U^{\mu}}{\partial\hat{n}}\,ds=\mu(\operatorname{int}(\gamma))=1.
\]
\end{enumerate}

\begin{theorem}[Gauss' Theorem; Theorem 2.1.1 of page 66 of \cite{SaffTotik2024}]\label{thm:gauss}
Let $R$ be a bounded open set with piecewise $C^1$ boundary, and $\mu$ a signed measure with compact support that is disjoint from $\partial R$. Then
\[
\mu(R)=\frac{1}{2\pi}\int_{\partial R}\frac{\partial U^{\mu}}{\partial\hat{n}}\,ds.
\]
Here $\frac{\partial}{\partial\hat{n}}$ denotes differentiation in the direction of the inner normal of $R$ and $ds$ indicates integration with respect to the arc length on $\partial R$.
\end{theorem}

Since $\gamma$ is an arbitrary analytic Jordan curve encircling $K$, by exhausting $\overline{\mathbb{C}}\setminus K$ (which is a simply connected domain in the Riemann sphere) by analytic Jordan domains $J_n$, with boundary curves $\gamma_n\subset\mathbb{D}\setminus K$ such that $\operatorname{int}(\gamma_n)\searrow K$, we find that $\nu(K)=\lim_{n}\nu(\operatorname{int}(\gamma_n))=1$. We have thus obtained our desired contradiction.

The main step that requires elaboration is Step 1 from which the rest of the proof follows in the way indicated by Steps 2--4. We complete the proof of Theorem~\ref{thm:main} by providing the details for Step 1 below.\label{step1-start}

Let
\[
\mathcal{K}(z,w):=\frac{1-|w|^2}{(1-\overline{w}z)(z-w)}
\]
for $z,w\in\mathbb{D}$. Then writing for a Borel probability measure $\rho$ of compact support contained in the unit disk
\[
C_{\mathbb{D}}^{\rho}(z):=\int_{\mathbb{D}}\mathcal{K}(z,w)\,d\rho(w),
\]
one has the relation
\[
C_{\mathbb{D}}^{\rho}(z)=C^{\rho}(z)+\int_{\mathbb{D}}\frac{\overline{w}}{1-\overline{w}z}\,d\rho(w)=C^{\rho}(z)-C^{\widehat{\rho}}(z),
\]
where for a finite Borel measure $\rho$ with compact support in the complex plane $\mathbb{C}$,
\[
C^{\rho}(z):=\int_{\mathbb{C}}\frac{1}{z-w}\,d\rho(w)
\]
is the Cauchy transform of $\rho$, and for a finite Borel measure $\rho$ with compact support in the unit disk $\mathbb{D}$, as before, $\widehat{\rho}$ denotes the Balayage of $\rho$ onto $\mathbb{T}^1$. It is well known that if $\rho_n$, $n=1,2,3,\dots$ is a sequence of finite compactly supported (in the complex plane $\mathbb{C}$) Borel measures, converging in the weak$^*$-topology of Borel measures on $\overline{\mathbb{C}}$ to the finite compactly supported Borel measure $\rho$, then $C^{\rho_n}\to C^{\rho}$ holds in $L_{\mathrm{loc}}^{1}(\mathbb{C})$ (for a proof of this fact, see e.g.\ the computation on page 2412 of \cite{Totik2019}). That is, for any fixed compact disk $\Delta(0,R):=\{z:|z|\le R\}$ one has
\[
\int_{\Delta(0,R)}|C^{\rho_n}(z)-C^{\rho}(z)|\,dA(z)\rightarrow_{n\rightarrow\infty} 0,
\]
where $dA$ is the two-dimensional Lebesgue measure. Consequently, if $\rho_n$, $n=1,2,\dots$ is a sequence of finite Borel measures, each with compact support in the unit disk $\mathbb{D}$, converging in the weak$^*$-topology of Borel measures on $\overline{\mathbb{C}}$ to the finite Borel measure $\rho$, with compact support in the unit disk $\mathbb{D}$, then one has that
\[
\int_{\Delta(0,r)}|C_{\mathbb{D}}^{\rho_n}(z)-C_{\mathbb{D}}^{\rho}(z)|\,dA(z)\rightarrow_{n\rightarrow\infty} 0
\]
for any $0<r<1$ fixed. Indeed, this follows from the same property for Cauchy transforms and the inequality
\begin{align*}
\int_{\Delta(0,r)}|C_{\mathbb{D}}^{\rho_n}(z)-C_{\mathbb{D}}^{\rho}(z)|\,dA(z)
&\le \int_{\Delta(0,r)}|C^{\rho_n}(z)-C^{\rho}(z)|\,dA(z)\\
&\quad +\int_{\Delta(0,r)}\left|\int_{\mathbb{D}}\frac{\overline{w}}{1-\overline{w}z}\,d\rho_n(w)-\int_{\mathbb{D}}\frac{\overline{w}}{1-\overline{w}z}\,d\rho(w)\right|dA(z).
\end{align*}
Noting that the term
\[
\int_{\Delta(0,r)}\left|\int_{\mathbb{D}}\frac{\overline{w}}{1-\overline{w}z}\,d\rho_n(w)-\int_{\mathbb{D}}\frac{\overline{w}}{1-\overline{w}z}\,d\rho(w)\right|dA(z)=o(1)
\]
for any fixed $0<r<1$ by the bounded convergence theorem and the weak$^*$-convergence $\rho_n\to\rho$. In other words,
\[
C_{\mathbb{D}}^{\rho_n}\rightarrow_{n\rightarrow\infty}C_{\mathbb{D}}^{\rho}
\]
holds in $L_{\mathrm{loc}}^{1}(\mathbb{D})$.

Since
\[
\frac{B_n'(z)}{nB_n(z)}=C_{\mathbb{D}}^{\mu_n}(z),
\]
and $\mu_n$ is a sequence of finite Borel measures, each compactly supported in the unit disk $\mathbb{D}$ and converging in the weak$^*$-topology of Borel measures on $\overline{\mathbb{C}}$ to $\mu$, which is compactly supported in the unit disk (in particular, $\operatorname{supp}(\mu)\subset K\subset\mathbb{D}$), we find that
\[
C_{\mathbb{D}}^{\mu_n}\rightarrow_{n\rightarrow\infty}C_{\mathbb{D}}^{\mu}
\]
in $L_{\mathrm{loc}}^{1}(\mathbb{D})$ and thus (after passing to an appropriate further subsequence), also (Lebesgue) almost everywhere in $\mathbb{D}$. It will thus be sufficient to show that the function $C_{\mathbb{D}}^{\mu}$ (which is holomorphic in $\mathbb{D}\setminus K$) has no zeros in $\mathbb{D}\setminus K$.

To this end, let $\kappa_n:=\sum_{j=1}^{m_n}\frac{p_{n,j}}{q_n}\delta_{z_{n,j}}$, $z_{n,j}\in K$, $j=1,\dots,m_n$ be a sequence of discrete probability measures for which $p_{n,j}\in\mathbb{N}_{\ge 1}$ for each $j=1,\dots,m_n$ and $q_n\in\mathbb{N}_{\ge 1}$ for each $n=1,2,3,\dots$, which converge to $\mu$ in the weak$^*$-topology of Borel measures on the Riemann sphere $\overline{\mathbb{C}}$. In other words, we take an approximation of $\mu$ by a sequence of rational (convex) linear combinations of delta measures, supported at points of $K$. Then for any point $z\in\mathbb{D}\setminus K$, we have
\[
C_{\mathbb{D}}^{\kappa_n}(z)=\int_{\mathbb{D}}\mathcal{K}(z,w)\,d\kappa_n(w)=\frac{1}{q_n}\frac{M_n'(z)}{M_n(z)},
\]
where
\[
M_n(z)=\prod_{j=1}^{m_n}\left(\frac{z-z_{n,j}}{1-\overline{z_{n,j}}z}\right)^{p_{n,j}}.
\]
By Theorem~\ref{thm:walsh}, $C_{\mathbb{D}}^{\kappa_n}(z)$ is a sequence of non-vanishing holomorphic functions in $\mathbb{D}\setminus K$, which (by a normal family argument and the weak$^*$-convergence) converges locally uniformly in $\mathbb{D}\setminus K$ to the holomorphic function $C_{\mathbb{D}}^{\mu}(z)\upharpoonright \mathbb{D}\setminus K$. By Hurwitz's theorem, we conclude that $C_{\mathbb{D}}^{\mu}(z)\upharpoonright \mathbb{D}\setminus K$ is non-vanishing in $\mathbb{D}\setminus K$ since it is clearly not identically zero there. Indeed,
\[
C_{\mathbb{D}}^{\mu}(z)=-2\frac{\partial}{\partial z}U_{\mathbb{D}}^{\mu}(z),
\]
where
\[
\frac{\partial}{\partial z}=\frac{1}{2}\left(\frac{\partial}{\partial x}-i\frac{\partial}{\partial y}\right),
\]
hence if $C_{\mathbb{D}}^{\mu}$ were identically zero in $\mathbb{D}\setminus K$, then since $-2\partial_z U_{\mathbb{D}}^{\mu}(z) = C_{\mathbb{D}}^{\mu}(z) = 0 = \overline{C_{\mathbb{D}}^{\mu}(z)} = -2\partial_{\overline{z}} U_{\mathbb{D}}^{\mu}(z)$, the Green potential $U_{\mathbb{D}}^{\mu}$ of $\mu$ would be constant in $\mathbb{D}\setminus K$, and since $U_{\mathbb{D}}^{\mu}$ has boundary values $0$ on the unit circle, and extends continuously there, this would imply that $U_{\mathbb{D}}^{\mu}$ is identically zero on $\mathbb{D}\setminus K$. But since $U_{\mathbb{D}}^{\mu}$ is nonnegative and superharmonic in the unit disk, the minimum principle implies that it must vanish identically all throughout the unit disk. Taking the distributional Laplacian, we find that $\mu\equiv 0$, contrary to fact ($\mu$ is a probability measure).

It follows that, after passing to an appropriate further subsequence,
\[
\frac{B_n'(z)}{nB_n(z)}
\]
converges (Lebesgue) almost everywhere in $\mathbb{D}\setminus K$ to the holomorphic function $C_{\mathbb{D}}^{\mu}(z)\upharpoonright \mathbb{D}\setminus K$, which is nowhere vanishing in $\mathbb{D}\setminus K$. Consequently,
\[
\frac{1}{n}\log\left|\frac{B_n'(z)}{nB_n(z)}\right|\rightarrow_{n\rightarrow\infty} 0
\]
holds for (Lebesgue) almost every $z\in\mathbb{D}\setminus K$, and we thus obtain the desired conclusion for Step 1. This completes the proof of Theorem~\ref{thm:main}.\label{step1-end}
\end{proof}

\section*{Statements and Declarations}
The author declares there is no conflict of interest regarding the publication of this paper.

\end{document}